\documentclass[12pt,a4paper]{article}
\usepackage[english]{}
\usepackage{amsfonts}
\usepackage{amsthm}
\usepackage[ansinew]{inputenc}
\usepackage{eufrak}
\usepackage[sumlimits,intlimits, namelimits]{amsmath}
\usepackage{graphicx}
\usepackage{flafter}
\usepackage{color}
\usepackage{multicol}
\usepackage{fancybox}
\usepackage{rotating}
\usepackage{amssymb}

\def\N{{\mathbb N}}
\def\Z{{\mathbb Z}}
\def\R{{\mathbb R}}

\def\weakstar{\stackrel{\ast}{\rightharpoonup}}

\def\mcL{\mathcal{L}}

\def\mutildeN{\tilde{\mu}_N}

\newcommand{\Proof}{\par\noindent{\bf Proof. }}
\newcommand{\eop}{\nopagebreak\hspace*{\fill}$\Box$\smallskip}

\newtheorem{theorem}{Theorem}[section]
\newtheorem{proposition}[theorem]{Proposition}

\newtheorem{lemma}[theorem]{Lemma}

\theoremstyle{definition}

\definecolor{grau}{gray}{0.80}

\title{Ground states of the 2D sticky disc model: fine properties and $N^{3/4}$ law for the deviation from the asymptotic Wulff shape}

\author{Bernd Schmidt\footnote{Universit{\"a}t Augsburg, Institut f{\"u}r Mathematik, 
Universit{\"a}tsstr.\ 14, 86159 Augsburg, Germany. {\tt bernd.schmidt@math.uni-augsburg.de}}}

\date{\today}

\begin{document}

\maketitle

\begin{abstract}
We investigate ground state configurations for a general finite number $N$ of particles of the Heitmann-Radin sticky disc pair potential model in two dimensions. Exact energy minimizers are shown to exhibit large microscopic fluctuations about the asymptotic Wulff shape which is a regular hexagon: There are arbitrarily large $N$ with ground state configurations deviating from the nearest regular hexagon by a number of $\sim N^{3/4}$ particles. We also prove that for any $N$ and any ground state configuration this deviation is bounded above by $\sim N^{3/4}$. As a consequence we obtain an exact scaling law for the fluctuations about the asymptotic Wulff shape. In particular, our results give a sharp rate of convergence to the limiting Wulff shape. 
\end{abstract}

\section{Introduction}

It is a fundamental problem in statistical and solid mechanics to derive material properties from atomistic interaction models. One of the important open questions is to explain theoretically why atoms at low temperature arrange on periodic lattices. In fact, at very low temperature, a system of atoms is not only observed to arrange on the atomic length-scale. The atoms are moreover seen to form large clusters whose overall polyhedral shape is given by the Wulff shape which is obtained through a surface energy minimization problem of an effective continuum model. 

The first crystallization result was obtained by Heitmann and Radin in \cite{HeitmannRadin:80}: For a very special pair interaction model in two dimensions, see \eqref{energy} and \eqref{radinpotential} below, they showed that ground state configurations at zero temperature are subsets of a triangular lattice with nearest-neighbour distance $1$. (See also \cite{Harborth} for an equivalent problem in a different context.) The overall shape of energy minimizers for this model was analyzed in \cite{AFS}: We proved that in the large $N$ limit, the unique rescaled asymptotic shape of ground states is given by a regular hexagon, which forms the Wulff shape of a limiting continuum variational problem. 

In this paper we address the problem of characterizing the discrete energy minimizers for finite $N$. In particular, we determine an exact scaling law for the microscopic fluctuations about the limiting Wulff shape: Ground state configurations deviate from the nearest regular hexagonal shape by at most $\sim N^{3/4}$ particles. This bound is sharp as there are highly degenerate discrete ground states for infinitely many $N$, which deviate form the nearest regular hexagon by $\sim N^{3/4}$ particles. It thus turns out that there are surprisingly large fluctuations even at zero temperature, as rearrangements of surface atoms would only lead to deviations of order $\sim N^{1/2}$. 

\begin{figure}[h!]
\begin{center}
\setlength{\unitlength}{1cm}
\begin{picture}(5,5)

\put(0.4,1.5){\vector(1,0){4.1}}
\put(0.5,0){\vector(0,1){4.5}}

\put(4.6,1.4){\footnotesize $r$}
\put(0.1,4.7){\footnotesize $V_{\rm HR}(r)$}

\put(0.2,1.4){\footnotesize $0$}
\put(0,3.9){\footnotesize $\infty$}
\put(0.4,4){\line(1,0){2.05}}
\put(-0.1,0.4){\footnotesize $-1$}
\put(0.4,0.5){\line(1,0){0.1}}

\put(2.52,1.6){\footnotesize $1$}

\put(2.5,4){\circle{0.1}}
\put(2.5,0.5){\circle*{0.1}}
\put(2.5,1.5){\circle{0.1}}

\multiput(0.5,4)(0.01,0){196}{\circle*{0.01}}
\multiput(2.55,1.5)(0.01,0){196}{\circle*{0.01}}

\multiput(2.45,4)(0,-0.10){35}{\circle*{0.01}}
\multiput(2.55,0.5)(0,0.10){10}{\circle*{0.01}}

\end{picture}
\end{center}
\caption{\label{F:potentials} The Heitmann-Radin sticky disc potential.}
\end{figure}
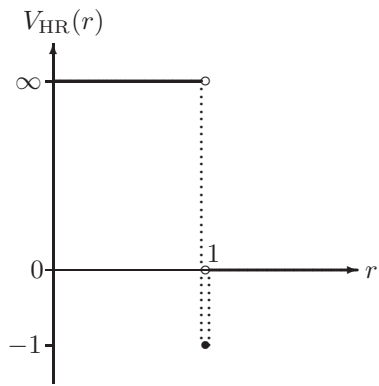

Our model energy, the Heitmann-Radin sticky disc model, assigns to $N$ particles in two dimensions with positions $x_1,\ldots,x_N\in\R^2$ the
pair potential energy 
\begin{align}\label{energy}
  E(x_1,\ldots,x_N) = \sum_{i\neq j} V_{\rm HR}(|x_i-x_j|),
\end{align}
with interaction potential given by
\begin{align}\label{radinpotential}
  V_{\rm HR}(r) = \left\{\begin{array}{ll}   +\infty, & 0\le r < 1, \\
                                       -1,      & r=1, \\
                                       0,       & r>1.\end{array}\right.
\end{align}
Such a potential can by interpreted as a very rough approximation of a hard core Lennard-Jones-type potential with infinitely short interaction length, see Figure \ref{F:potentials}.

Based on the crystallization result in \cite{Harborth, HeitmannRadin:80}, for this model it was shown \cite{AFS} that ground states approach a unique hexagonal shape, in the sense that the re-scaled empirical measure of any $N$-particle minimizer
$(x_1^{(N)},\ldots,x_N^{(N)})$, 
\begin{align} \label{C:radonmeasure}
  \mu_N = N^{-1} \sum_{j=1}^N \delta_{N^{-1/2}x_j^{(N)}},
\end{align}
satisfies, up to translation and rotation,
\begin{align} \label{limitmeasure}
  \mu_N \weakstar \rho_0 \chi_{h},
\end{align}
where $\chi_{h}$ is the characteristic function of the regular hexagon $h$ with corners $\pm \frac{1}{\sqrt{3}}{\footnotesize 1 \choose 0}$, $\pm \frac{1}{2\sqrt{3}} {\footnotesize 1 \choose \sqrt{3}}$, $\pm \frac{1}{2\sqrt{3}} {\footnotesize -1 \choose \sqrt{3}}$, and $\rho_0 = \frac{2}{\sqrt{3}}$ is the density of atoms in a triangular lattice with nearest-neighbour distance $1$. 

For finite $N$ we will measure deviations $\|\mu_N - \rho_0 \chi_h\|$ from the nearest regular hexagon $h$ in terms of a suitable norm, see \eqref{eq:flat-norm}. Our main results are as follows.

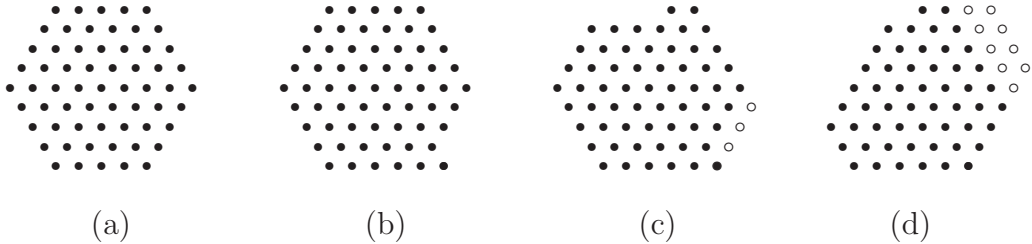
\begin{figure}[h!]
\begin{center}
\setlength{\unitlength}{0.3cm}
\begin{picture}(48,10.5)

\put(0,3){%
\setlength{\unitlength}{0.3cm}
\begin{picture}(10,8)
\multiput(3,0)(1,0){5}{\circle*{0.30}}
\multiput(2.5,0.866)(1,0){6}{\circle*{0.30}}
\multiput(2,1.732)(1,0){7}{\circle*{0.30}}
\multiput(1.5,2.598)(1,0){8}{\circle*{0.30}}
\multiput(1,3.464)(1,0){9}{\circle*{0.30}}
\multiput(1.5,4.33)(1,0){8}{\circle*{0.30}}
\multiput(2,5.196)(1,0){7}{\circle*{0.30}}
\multiput(2.5,6.062)(1,0){6}{\circle*{0.30}}
\multiput(3,6.928)(1,0){5}{\circle*{0.30}}

\end{picture}
}

\put(12,3){%
\setlength{\unitlength}{0.3cm}
\begin{picture}(11,8)
\multiput(3,0)(1,0){6}{\circle*{0.30}}
\multiput(2.5,0.866)(1,0){6}{\circle*{0.30}}
\multiput(2,1.732)(1,0){7}{\circle*{0.30}}
\multiput(1.5,2.598)(1,0){8}{\circle*{0.30}}
\multiput(1,3.464)(1,0){9}{\circle*{0.30}}
\multiput(1.5,4.33)(1,0){8}{\circle*{0.30}}
\multiput(2,5.196)(1,0){7}{\circle*{0.30}}
\multiput(2.5,6.062)(1,0){6}{\circle*{0.30}}
\multiput(3,6.928)(1,0){5}{\circle*{0.30}}

\multiput(8,0)(0.5,0.866){1}{\circle*{0.30}}
\end{picture}
}

\put(24,3){%
\setlength{\unitlength}{0.3cm}
\begin{picture}(11,8)
\multiput(3,0)(1,0){6}{\circle*{0.30}}
\multiput(2.5,0.866)(1,0){6}{\circle*{0.30}}
\multiput(2,1.732)(1,0){7}{\circle*{0.30}}
\multiput(1.5,2.598)(1,0){8}{\circle*{0.30}}
\multiput(1,3.464)(1,0){9}{\circle*{0.30}}
\multiput(1.5,4.33)(1,0){8}{\circle*{0.30}}
\multiput(2,5.196)(1,0){7}{\circle*{0.30}}
\multiput(2.5,6.062)(1,0){6}{\circle*{0.30}}
\multiput(6,6.928)(1,0){2}{\circle*{0.30}}

\multiput(8,0)(0.5,0.866){4}{\circle{0.30}}
\end{picture}
}

\put(35,3){%
\setlength{\unitlength}{0.3cm}
\begin{picture}(11,8)
\multiput(3,0)(1,0){6}{\circle*{0.30}}
\multiput(2.5,0.866)(1,0){6}{\circle*{0.30}}
\multiput(2,1.732)(1,0){7}{\circle*{0.30}}
\multiput(2.5,2.598)(1,0){7}{\circle*{0.30}}
\multiput(3,3.464)(1,0){7}{\circle*{0.30}} \multiput(10,3.464)(1,0){1}{\circle{0.30}}
\multiput(3.5,4.33)(1,0){6}{\circle*{0.30}} \multiput(9.5,4.33)(1,0){2}{\circle{0.30}}
\multiput(4,5.196)(1,0){5}{\circle*{0.30}} \multiput(9,5.196)(1,0){2}{\circle{0.30}}
\multiput(4.5,6.062)(1,0){4}{\circle*{0.30}} \multiput(8.5,6.062)(1,0){2}{\circle{0.30}}
\multiput(6,6.928)(1,0){2}{\circle*{0.30}} \multiput(8,6.928)(1,0){2}{\circle{0.30}}

\multiput(8,0)(0.5,0.866){4}{\circle*{0.30}}
\end{picture}
}

\put(5,0){(a)}
\put(17,0){(b)}
\put(29,0){(c)}
\put(40,0){(d)}

\end{picture}
\end{center}
\caption{\label{F:examples} (a) The unique (up to rotations and translations) ground state for 61 atoms, cf.\ Proposition \ref{prop:naturalcandidate}. (b) A ground state for 62 atoms, cf.\ Proposition \ref{prop:unique-gs}. (c) Another ground state for 62 atoms obtained from (b) by rearrangement of surface atoms. New positions are indicated by white lattice sites. (d) Yet another ground state for 62 atoms obtained from (c) by rearrangement of a large surface layer of atoms.} 
\end{figure}

\begin{itemize}
\item A simple observation shows that there always exists a ground state configuration which agrees with a suitably chosen regular hexagon up to a surface contribution of $\sim N^{1/2}$ atoms, i.e., $\|\mu_N - \rho_0 \chi_h\| \le N^{-1/2}$, see Figures \ref{F:examples}(a),(b),(c) and \ref{F:hex}. 

\item For many arbitrarily large $N$ there exist ground states which are highly degenerate: they differ
from the nearest regular hexagon by $\sim N^{3/4}$ atoms, i.e., $\|\mu_N - \rho_0 \chi_h\| \ge c N^{-1/4}$, see Figure \ref{F:examples}(d) (and Theorem \ref{theo:lowerbound}). 

\item Finally, this rate is sharp, in the sense that $\| \mu_N - \rho_0 \chi_h\| \le C N^{-1/4}$ for all $N$ and any ground state (see Theorem \ref{theo:rateupper}). 
\end{itemize}

It is interesting to note that there are also infinitely many $N$ for which the ground state is (up to rotation and translation) unique and for which the ground states of $N+1$ particles are highly degenerate. 

As an offspring, we obtain an alternative, perhaps conceptually less transparent but more elementary, proof of our result in \cite{AFS} of convergence to the hexagonal Wulff shape, which also provides a sharp rate of convergence. 

We also remark that obviously our results hold true for any pair interaction model with interaction range less than the next nearest neighbor distance $\sqrt{3}$, for which ground states are known to crystallize on the triangular lattice. This is e.g.\ the case for the Radin soft disc potential $V_{\rm R}$, for which $V_{\rm R}(r) = 24r - 25$ on $[1, \frac{25}{24}]$ and which equals $V_{\rm HR}$ elsewhere, see \cite{Radin:81}, and thus for any $V$ with $V_{\rm R} \le V \le V_{\rm HR}$. However, it seems that there are no examples known which allow for an elastic range, i.e., $V$ quadratic in a neighborhood of the equilibrium point $1$, as e.g.\ the short-range soft potentials considered in \cite[Section 2]{AFS}, although we conjecture our results also to remain valid for these models. In fact, we expect that even long-range Lennard-Jones type interactions lead to Wulff shapes in the infinite particle limit. Theil remarkably succeded in proving crystallization under crystalline or periodic boundary conditions for a class of such models \cite{Theil}. In our context, however, such a result would have to be proven, at least asymptotically, for the free boundary value problem, which appears to be a very difficult problem. 

Our results explain, in a model problem, fluctuations of atomic configurations at zero temperature. The rapid change of unique and highly degenerate ground states as $N$ varies leads to an interesting discrepancy between the canonical and grandcanonical Gibbs ensemble restricted to absolute minimizers of the potential energy. It would be interesting to explore to what extend these results persist at low but finite temperature. In this regime the canonical Gibbs ensemble is expected to behave more regularly since even for particle numbers $N$ with unique ground states a construction as in the proof of Theorem \ref{theo:lowerbound} shows that the second lowest energy level is highly degenerate. In line with our analysis at zero temperature, we expect to observe surprisingly large entropic effects for the statistical fluctuations about the hexagonal Wulff also at finite temperature. 

\section{Analysis of ground states}

Let ${\cal S}_N=\{x_1,\ldots,x_N\}$ denote a general configuration of atomic positions $x_i\in\R^2$ with finite energy. Hence $|x_i-x_j|\ge 1$ for all $i\neq j$.   
By viewing the atomic positions as vertices that are linked to each other by an atomic bond whenever their distance is equal to $1$, we have assigned a graph to every configuration of atoms, whose set of vertices is ${\cal S}_N$ and whose set of edges ${\cal B}$ consists of their atomic bonds.  

An element of ${\cal B}$ is called a boundary edge and said to belong to $\partial {\cal B}$ if there are less 
than two vertices that are connected with both of its endpoints. We call those $x \in {\cal S}_N$ for which $\# \{ y \in {\cal S}_N : |y - x| = 1 \} \le 5$ boundary points and denote the set of boundary points by $\partial {\cal S}_N$.

The main result of Heitmann and Radin in \cite{HeitmannRadin:80} states that ground states are subsets of a 
triangular lattice with minimal interatomic distance $1$. Of course any rigid motion maps ground states into ground states, and we will therefore without loss of generality only consider atomic configurations ${\cal S}_N$ which are subsets of the lattice 
\begin{align}\label{lattice}
  \mcL := \{ m e_1 + n e_2 : m, \, n\in\Z\}, \;\;
  e_1=\begin{pmatrix} 1 \\ 0 \end{pmatrix}, \;\;
  e_2=\mbox{$\frac12$}\begin{pmatrix} 1 \\ \sqrt{3} \end{pmatrix}.
\end{align}
For later reference we include a precise statement of the Heitmann-Radin results.

\begin{theorem}[\cite{Harborth,HeitmannRadin:80}]\label{theo:HeitmannRadin}
Let ${\cal S}_N$ be a minimal energy configuration of $N \ge 3$ atoms with set of atomic bonds ${\cal B}$.
\begin{itemize}
\item[(i)] Then $\# {\cal B} = \lfloor 3 N - \sqrt{12 N - 3} \rfloor$.
\item[(ii)] Up to a rigid displacement the elements of $\partial{\cal B}$ form a simple closed polygon $P$ 
with vertices on ${\cal L}$. ${\cal S}_N$ consists precisely of all the points inside and on $P$. 
\item[(iii)] $\# \partial {\cal S}_N = \# \partial{\cal B} = - \lfloor 3 - \sqrt{12 N - 3} \rfloor$.
\end{itemize}
\end{theorem}

As an immediate consequence of (i) we obtain that the ground state energy for $N$ atoms is given by 
\begin{align}\label{eq:gsenergy} 
  -2 \# {\cal B} 
  = -2 \left\lfloor 3N - \sqrt{12N - 3} \right\rfloor 
  = - 6N + 2 \left\lceil \sqrt{12N - 3} \right\rceil. 
\end{align}

We are interested in the ground states for large finite $N$. So in the following results it is no loss of generality to (tacitly) assume that $N$ be sufficiently large.

Let ${\cal S}_N$ be a ground state. We choose an ordering of the elements of $\partial {\cal S}_N = {\cal S}_N \cap P$ by 
running through $P$ counter-clockwise as $v_1, v_2, \ldots, v_m, v_{m+1} = v_1$. For $v \in \partial {\cal S}_N$ we 
denote by $\varphi_v$ the interior angle of $P$ at $v$ (taking values in $\{\frac{\pi}{3}, \frac{2\pi}{3}, \pi, 
\frac{4\pi}{3} \}$) and set $\varphi_i = \varphi_{v_i} = <\!\!\!) (v_{i + 1} - v_i, v_{i-1} - v_i)$. 

\begin{lemma}\label{lemma:structurei}
Let ${\cal S}_N$ be a ground state. 
\begin{itemize}
\item[(i)] Then $\varphi_i = \frac{4\pi}{3}, \varphi_{i+1} = \ldots = \varphi_{i+j} = \pi$ implies that 
$\varphi_{i+j+1} \ne \frac{4\pi}{3}$. (Here and in the sequel all indices are understood modulo $m$.) 
\item[(ii)] If $\varphi_i = \frac{\pi}{3}$, then $\varphi_j \ne \frac{4\pi}{3}$ for $j \notin \{i-1,i+1\}$. 
\item[(iii)] There is at most one $i \in \{1, \ldots, m\}$ such that $\varphi_i = \frac{\pi}{3}$.
\item[(iv)] If $\varphi_i = \frac{\pi}{3}$, then $\varphi_{i-1}, \varphi_{i+1} \in \{ \pi, \frac{4\pi}{3} \}$ 
and not both of these angles are equal to $\pi$.
\end{itemize}
\end{lemma}

\Proof
Note first that for any configuration of atoms filling the lattice inside and on some polygon $P$ whose segments are aligned with the nearest neighbor bonds in ${\cal L}$, there always exists a boundary vertex $v$ of degree ${\rm deg}(v) \le 3$ which lies in a prescribed half plane that intersects $P$. (Simply take a suitable corner of $P$.) 

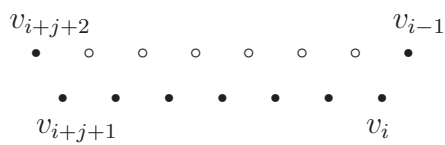
\begin{figure}[h!]
\begin{center}
\setlength{\unitlength}{0.7cm}
\begin{picture}(8,3)
\put(0.5,1.72){\circle*{0.17}}
\multiput(1,0.86)(1,0){7}{\circle*{0.17}}
\multiput(1.5,1.72)(1,0){6}{\circle{0.17}}
\put(7.5,1.72){\circle*{0.17}}
\put(0.5,0.2){$v_{i+j+1}$}
\put(6.7,0.2){$v_{i}$}
\put(0, 2.2){$v_{i+j+2}$}
\put(7.2, 2.2){$v_{i-1}$}
\end{picture}
\end{center}
\caption{\label{fig:freeslots} Free lattice slots between atom $v_{i-1}$ and $v_{i+j+2}$.}
\end{figure}

(i) If this were not the case, i.e., we had a configuration as sketched in Figure \ref{fig:freeslots},

\noindent we could successively move atoms of degree at most $3$ in between the atoms $v_{i-1}$ and $v_{i+j+2}$. In each step 
we would loose at most three atomic bonds by removing the atom but---except for the last move---also gain three new bonds by 
putting the atom back in. In the last step, however we would loose three and gain four new bonds, thus 
contradicting maximality of the number of bonds. (Since $j \ll N$, we do not have to move the atoms $v_{i-1}, 
\ldots, v_{i+j+2}$ themselves.)

(ii) If $\varphi_i = \frac{\pi}{3}$ and $\varphi_j = \frac{4\pi}{3}$ for some $j \notin \{i-1,i+1\}$, then 
moving the atom $v_i$ to the point of ${\cal L}$ adjacent to $v_{j-1}$, $v_j$ and $v_{j+1}$ 
gives three new bonds while only destroying two old ones, in contradiction to maximality. 

(iii) If $\varphi_i = \varphi_j = \frac{\pi}{3}$ for $i \ne j$, by (ii) we have $\varphi_k \in \{ 
\frac{\pi}{3}, \frac{2\pi}{3}, \pi \}$ for all $k \notin \{i-1,i+1\} \cup \{j-1,j+1\}$. This implies that for at most 
six of these $k$ the angle $\varphi_k$ can be less than $\pi$, because $<\!\!\!) (v_{i+2} - v_{i+1}, v_{i-2} - v_{i-1}), <\!\!\!) (v_{j+2} - v_{j+1}, v_{j-2} - v_{j-1}) \le \pi$ and $P$ forms a simple closed polygon. We therefore find some $k_0 \notin \{i-1,i,i+1\} \cup 
\{j-1,j,j+1\}$ with $\varphi_{k_0} = \pi$. Now we can move the atoms $v_i$ and $v_j$ to the two points of ${\cal 
L} \setminus {\cal S}_N$ adjacent to $v_{k_0}$ with a net gain of one atomic bond. 

(iv) It is easy to see that $\varphi_{i-1}, \varphi_{i+1} \ne \frac{\pi}{3}$. An argument as in the proof of 
(iii) shows that there is a $k$ such that $\varphi_k = \varphi_{k+1} = \pi$. If $\varphi_{i-1} = \varphi_{i+1} = 
\pi$, then $v_{i-1}, v_i, v_{i+1}$ could be moved to the three points of ${\cal L} \setminus {\cal S}_N$ that are 
adjacent to $v_{k}$ or $v_{k+1}$. This would break seven atomic bonds but create eight new ones. Now 
if without loss of generality $\varphi_{i+1}$ is assumed to be $\frac{2\pi}{3}$, we similarly obtain $5 - 4 = 1$ 
new bond by moving of $v_i$ and $v_{i+1}$ to the two points of ${\cal L} \setminus {\cal S}_N$ that are adjacent 
to $v_{k}$.
\eop

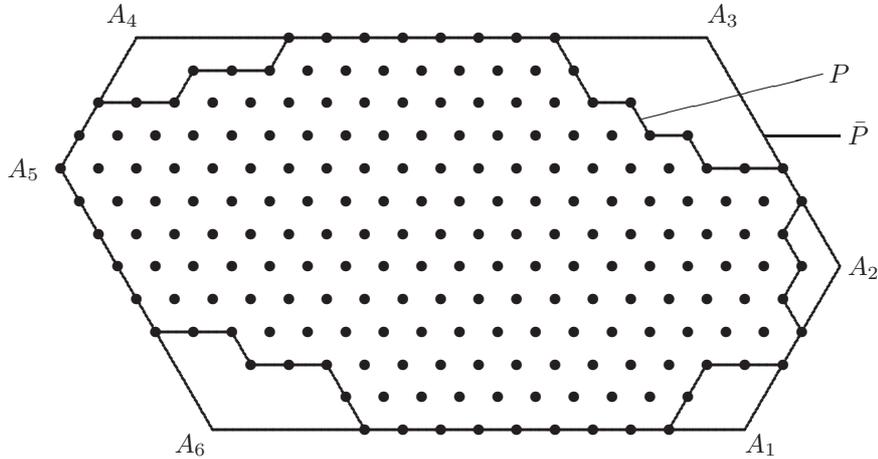
\begin{figure}[h!]
\begin{center}
\setlength{\unitlength}{1cm}
\begin{picture}(12,6.1)
\put(0,-0.7){%
\setlength{\unitlength}{1cm}
\begin{picture}(12,6.8)

\multiput(3,1)(0.01,0){701}{\circle*{0.01}}
\multiput(3,1)(-0.01,0.0173){201}{\circle*{0.01}}
\multiput(1,4.464)(0.01,0.0173){101}{\circle*{0.01}}
\multiput(2,6.196)(0.01,0){751}{\circle*{0.01}}
\multiput(10,1)(0.01,0.0173){126}{\circle*{0.01}}
\multiput(9.5,6.196)(0.01,-0.0173){176}{\circle*{0.01}}

\multiput(5,1)(0.5,0){9}{\circle*{0.15}}
\multiput(4.75,1.433)(0.5,0){10}{\circle*{0.15}}
\multiput(3.5,1.866)(0.5,0){15}{\circle*{0.15}}
\multiput(2.25,2.299)(0.5,0){18}{\circle*{0.15}}
\multiput(2,2.732)(0.5,0){18}{\circle*{0.15}}
\multiput(1.75,3.165)(0.5,0){19}{\circle*{0.15}}
\multiput(1.5,3.598)(0.5,0){19}{\circle*{0.15}}
\multiput(1.25,4.031)(0.5,0){20}{\circle*{0.15}}
\multiput(1,4.464)(0.5,0){20}{\circle*{0.15}}
\multiput(1.25,4.895)(0.5,0){17}{\circle*{0.15}}
\multiput(1.5,5.330)(0.5,0){15}{\circle*{0.15}}
\multiput(2.75,5.763)(0.5,0){11}{\circle*{0.15}}
\multiput(4,6.196)(0.5,0){8}{\circle*{0.15}}

\multiput(5,1)(-0.01,0.0173){50}{\circle*{0.01}}
\multiput(4.5,1.866)(-0.01,0){100}{\circle*{0.01}}
\multiput(3.5,1.866)(-0.01,0.0173){25}{\circle*{0.01}}
\multiput(3.25,2.299)(-0.01,0){100}{\circle*{0.01}}

\multiput(1.5,5.330)(0.01,0){100}{\circle*{0.01}}
\multiput(2.5,5.330)(0.01,0.0173){25}{\circle*{0.01}}
\multiput(2.75,5.763)(0.01,0){100}{\circle*{0.01}}
\multiput(3.75,5.763)(0.01,0.0173){25}{\circle*{0.01}}

\multiput(7.5,6.196)(0.01,-0.0173){50}{\circle*{0.01}}
\multiput(8,5.330)(0.01,0){50}{\circle*{0.01}}
\multiput(8.5,5.330)(0.01,-0.0173){25}{\circle*{0.01}}
\multiput(8.75,4.895)(0.01,0){50}{\circle*{0.01}}
\multiput(9.25,4.895)(0.01,-0.0173){25}{\circle*{0.01}}
\multiput(9.5,4.464)(0.01,0){100}{\circle*{0.01}}

\multiput(10.75,4.031)(-0.01,-0.0173){25}{\circle*{0.01}}
\multiput(10.5,3.598)(0.01,-0.0173){25}{\circle*{0.01}}
\multiput(10.75,3.165)(-0.01,-0.0173){25}{\circle*{0.01}}
\multiput(10.5,2.732)(0.01,-0.0173){25}{\circle*{0.01}}

\multiput(10.5,1.866)(-0.01,0){100}{\circle*{0.01}}
\multiput(9.5,1.866)(-0.01,-0.0173){50}{\circle*{0.01}}

\put(10,0.7){\footnotesize $A_1$}
\put(11.35,3.05){\footnotesize $A_2$}
\put(9.5,6.4){\footnotesize $A_3$}
\put(1.6,6.4){\footnotesize $A_4$}
\put(0.3,4.35){\footnotesize $A_5$}
\put(2.5,0.7){\footnotesize $A_6$}

\put(8.625,5.113){\line(4,1){2.4}}
\put(11.1,5.6){\footnotesize $P$}
\put(10.25,4.895){\line(1,0){1}}
\put(11.35,4.78){\footnotesize $\bar{P}$}
\end{picture}
}
\end{picture}
\end{center}
\caption{\label{fig:hexaround} Configuration with surrounding hexagon.}
\end{figure}

If there is an $i$ with $\varphi_i = \frac{\pi}{3}$, it follows from Lemma \ref{lemma:structurei} that deleting 
the atom $v_i$ we are left with a graph whose boundary polygon $\bar{P}$ is in fact a convex hexagon. In case 
there is no such $i$, we define $\bar{P}$ to be the smallest convex hexagon whose segments are aligned with the bonds in 
${\cal L}$ and that contains ${\cal S}_N$, i.e., the interior of $\bar{P}$ is the intersection of all half planes $\{x : 
x\cdot R_{\frac{(2i+1)\pi}{6}} {e}_1 \le c\}$, $i=0,\ldots,5$, $c\in\R$, containing ${\cal S}_N$, see Figure \ref{fig:hexaround}. (Here $R_{\theta}$ denotes the rotation about the angle $\theta$.) 

Let $A_1,\ldots,A_6 \in {\cal L}$ be the corners of $\bar{P}$ such that the segments 
$a_i := [A_i, A_{i+1}]$ are numbered such that $(A_{i+1} - A_i) \cdot R_{\frac{i\pi}{3}} {e}_1 = |A_{i+1} - A_i|$ (mod $6$).  Denote the points of ${\cal L}$ inside or on $\bar{P}$ by $\mathcal{T}$ and recall the connectedness definition from \cite{AFS}: a finite set ${\cal S}\subset\R^2$ is called connected if for any two $x,y\in{\cal S}$ there exist $x_0,\ldots,x_N\in{\cal S}$ such that $x_0=x$, $x_N=y$, and the distance between successive points $x_{j-1}$, $x_j$ lies within the interaction range of the potential, i.e. in our case $|x_j-x_{j-1}|=1$ for all $j$.   

\begin{lemma}\label{lemma:structureii}
Let ${\cal S}_N$ be a ground state with $\varphi_i \ne \frac{\pi}{3}$ for all $i$. 
\begin{itemize}
\item[(i)] There are points $P_i, P_i' \in a_i \cap {\cal L}$ such that the segments $(A_i, P_i)$ and $(P_i', A_{i+1})$ are disjoint and 
$$ {\cal S}_N \cap a_i = {\cal L} \cap [P_i, P_i']. $$ 
\item[(ii)] There are at most six connected components of $\mathcal{T} \setminus {\cal S}_N$, each of them 
containing precisely one corner $A_i$.
\item[(iii)] Label the set of atoms in these components accordingly with ${\cal C}_i$ and denote their respective boundary polygons by $C_i$, $i \in \{1, \ldots, 6\}$. Then $C_i$ is given by the two segments $[P_{i-1}', A_i], [A_i, P_i]$ and a `staircase' between $P_i$ and $P_{i-1}'$, i.e., a polygonal path whose segments are parallel to $a_{i-1}$ or $a_i$. 
\item[(iv)] $\# ({\cal C}_1 \cup \ldots \cup {\cal C}_6) \le |A_{i+1} - A_i|$ for all $i = 1, \ldots, 6$.
\end{itemize}
\end{lemma}

\Proof 
(i) Let $P_i$, resp.\ $P_i'$, be the point in $a_i \cap {\cal S}_N$ that is closest to $A_i$, resp.\ $A_{i+1}$. Then in fact ${\cal S}_N \cap a_i = {\cal L} \cap [P_i, P_i']$ for otherwise there existed lattice points $v_i, v_{i+j}$ on the polygonal path $P$ between $P_i$ and $P_i'$ with $\varphi_i = \frac{4\pi}{3}, \varphi_{i+1} = \ldots = \varphi_{i+j-1} = \pi, \varphi_{i+j} = \frac{4\pi}{3}$, contradicting Lemma \ref{lemma:structurei}(i).

(ii) This is a direct consequence of (i) and Theorem \ref{theo:HeitmannRadin}(ii). 

(iii) If ${\cal C}_i \ne \emptyset$, then $P_{i-1}' \ne A_i \ne P_i$. Let $P_{i-1}' = v_j, P_i = v_{j+r}$. Then $\varphi_j = \frac{2\pi}{3} = \varphi_{j+r}$ and for all $k = j, j+1, \ldots, j+r$ we have $\varphi_k \in \{\frac{2\pi}{3}, \pi, \frac{4\pi}{3} \}$. Let $k_l$ be the subsequence of those indices for which $\varphi_{k_l} \ne \pi$. Then the angles $\varphi_{k_l}$ alternate between $\frac{2\pi}{3}$ and $\frac{4\pi}{3}$, i.e., $\varphi_{k_l} = \frac{2\pi}{3} \iff \varphi_{k_{l+1}} = 
\frac{4\pi}{3}$, for otherwise we could find $l$ such that $\varphi_{k_l} = \frac{4\pi}{3} = \varphi_{k_{l+1}}$ in contradiction to Lemma \ref{lemma:structurei}(i). (Observe that $\#\{l : \varphi_{k_l} = \frac{2\pi}{3} \} = \#\{ l : \varphi_{k_l} = \frac{4\pi}{3} \} + 1$.) This implies the claim.

(iv) If $\# ({\cal C}_1 \cup \ldots \cup {\cal C}_6) \ge |A_{i+1} - A_i| + 1 = \# (a_i \cap {\cal L})$, then we 
can successively move the atoms of our given ground state that lie on the segment $a_i$ into the set $({\cal C}_1 \cup \ldots \cup {\cal C}_6) \setminus 
(a_i \cap {\cal L})$. This is possible since the number of those atoms is $\# ( (a_i \cap {\cal L}) \setminus ({\cal C}_1 \cup \ldots \cup {\cal C}_6) )$, which by assumption is not greater than  $ \# ( ({\cal C}_1 \cup 
\ldots \cup {\cal C}_6) \setminus (a_i \cap {\cal L}) )$. By (iii) these moves can be done in such a way 
that we gain three new atomic bonds in each step. Choosing an atom with the least neighbors among the atoms on 
$a_i$, we will also loose only three atomic bonds in every step, except for the last one. In the last step, however, i.e., if there is just one atom left to move, we will only loose two atomic bonds, thus contradicting that ${\cal S}_N$ was a ground state. 
\eop

The ground state will in general not be unique, in particular, the sets ${\cal C}_i$ introduced in the previous lemma may 
not consist of a single row of atoms. However, from any ground state ${\cal S}_N$ we can always construct a `normalized' ground state ${\cal S}_N^{\rm (n)}$ by moving as many atoms as possible from $\bar{P}$ to $({\cal C}_1 \cup \ldots \cup {\cal C}_6) \setminus \bar{P}$ as described in the proof of Lemma \ref{lemma:structureii}(iv). Moving the exterior atoms lying on $\bar{P}$, we may in addition assume that in fact ${\cal C}_2 = \ldots = {\cal C}_6 = \emptyset$ and $P = [A_1 + {e}_2, A_2] \cup [A_2, A_3] \cup \ldots \cup [A_5, A_6] \cup [A_6, P_6'] \cup [P_6', P_6' + {e}_2] \cup[P_6' + {e}_2, A_1 + {e}_2]$. 

We now quantify the deviation between the empirical measure of a finite-$N$ ground state and its asymptotic hexagonal shape (see \eqref{limitmeasure}). 

This deviation can be quantified e.g.\ in terms of the flat norm 
\begin{align}\label{eq:flat-norm}
  \|\mu\| 
  := \sup \left\{ \int \varphi \, d\mu : \mbox{ Lipschitz with } |\varphi| \le 1 
     \mbox{ and Lipschitz constant } \le 1 \right\}. 
\end{align}
Note that alternatively, instead of measuring the rate of convergence of the $\mu_N$ in terms the flat norm, we could start from a related absolute continuous measure
(as used for various purposes in \cite{AFS}) such as
$$
    \mutildeN = N^{-1} \sum_{x\in {\cal S}_N} \frac{\chi_{V_*(x)}}{\mbox{area}(V_*(x))},
$$
where $V_*(x)$ is the Voronoi cell of the point $x$ with respect to the re-scaled triangular lattice $N^{-1/2}{\cal L}$.  
We could then measure deviations to the limiting function $\mu$ in terms of the standard $L^1$-distance. Due to the fact that in a ground state the rescaled atomic positions lie on $\frac{1}{\sqrt{N}} {\cal L}$, the rate of convergence of $\mu_N$ to $\mu$ cannot be faster than $N^{-1}$. On the other hand we have 
$$ \left| \|\mu_N - \mu\| - \|\tilde{\mu}_N - \mu\|_{L^1} \right| \le C N^{-1}. $$ 
So both approaches will yield equivalent results. 

In passing from ${\cal S}_N$ to ${\cal S}_N^{\rm (n)}$ we have moved no more than $C N^{\frac{1}{2}}$ atoms, so 
\begin{align}\label{eq:normalizing}
  \| \mu_N - \mu_N^{\rm (n)} \| \le CN^{-\frac{1}{2}}, 
\end{align}
where $\mu_N^{\rm (n)}$ denotes the rescaled empirical measure of the normalized ground state configuration ${\cal S}_N^{\rm (n)}$. In fact, fluctuations of the order $N^{-\frac{1}{2}}$ are to be expected due to surface rearrangements. In the following we will see that the optimal rate of convergence is, somewhat surprisingly, only of order $N^{-\frac{1}{4}}$. We begin by deriving a corresponding upper bound.

\begin{theorem}\label{theo:rateupper} Suppose ${\cal S}_N \subset {\cal L}$ is a ground state. There is a constant $C$, independent of $N$, such that $\mu_N$, after a suitable translation, satisfies 
$$\| \mu_N - \mu \| \le C N^{-\frac{1}{4}}.$$
\end{theorem}

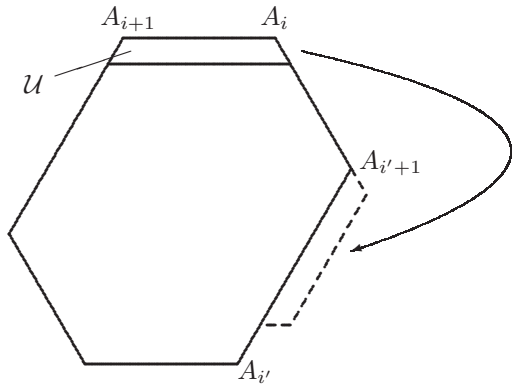
\begin{figure}[h!]
\begin{center}
\setlength{\unitlength}{1cm}
\begin{picture}(8,5.3)
\put(1,-0.7){%
\begin{picture}(8,6)

\multiput(1,1)(0.01,0){201}{\circle*{0.01}}
\multiput(1,1)(-0.01,0.0173){101}{\circle*{0.01}}
\multiput(0,2.732)(0.01,0.0173){151}{\circle*{0.01}}
\multiput(1.5,5.327)(0.01,0){201}{\circle*{0.01}}
\multiput(3.5,5.327)(0.01,-0.0173){101}{\circle*{0.01}}
\multiput(3,1)(0.01,0.0173){151}{\circle*{0.01}}

\multiput(1.3,4.981)(0.01,0){241}{\circle*{0.01}}

\multiput(4.5,3.595)(0.01,-0.0173){6}{\circle*{0.01}}
\multiput(4.6,3.422)(0.01,-0.0173){6}{\circle*{0.01}}

\multiput(4.7,3.249)(-0.01,-0.0173){6}{\circle*{0.01}}
\multiput(4.6,3.076)(-0.01,-0.0173){6}{\circle*{0.01}}
\multiput(4.5,2.903)(-0.01,-0.0173){6}{\circle*{0.01}}
\multiput(4.4,2.73)(-0.01,-0.0173){6}{\circle*{0.01}}
\multiput(4.3,2.557)(-0.01,-0.0173){6}{\circle*{0.01}}
\multiput(4.2,2.384)(-0.01,-0.0173){6}{\circle*{0.01}}
\multiput(4.1,2.211)(-0.01,-0.0173){6}{\circle*{0.01}}
\multiput(4.0,2.038)(-0.01,-0.0173){6}{\circle*{0.01}}
\multiput(3.9,1.865)(-0.01,-0.0173){6}{\circle*{0.01}}
\multiput(3.8,1.692)(-0.01,-0.0173){6}{\circle*{0.01}}

\multiput(3.7,1.517)(-0.01,0){12}{\circle*{0.01}}
\multiput(3.5,1.517)(-0.01,0){12}{\circle*{0.01}}

\put(4.6,3.595){\footnotesize $A_{i'+1}$}
\put(3,0.8){\footnotesize $A_{i'}$}
\put(1.2,5.5){\footnotesize $A_{i+1}$}
\put(3.3,5.5){\footnotesize $A_{i}$}
\put(1.6,5.15){\line(-3,-1){1}}
\put(0.2,4.6){\footnotesize ${\cal U}$}

\qbezier(3.85,5.15)(9,4)(4.5,2.5)
\put(4.5,2.5){\vector(-3,-1){.01}}
\end{picture}
}
\end{picture}
\end{center}
\caption{\label{fig:nolossii} Moving a large set of atoms keeping the number of bonds fixed.}
\end{figure}

\Proof By \eqref{eq:normalizing}, without loss of generality assume that ${\cal S}_N = {\cal S}_N^{\rm (n)}$. Suppose that $a_i$, $i \in \{1, \ldots, 6\}$, and $a_{i'}$, $i' \in \{1, \ldots, 5\}$ satisfy $|a_{i'}| - |a_i| \ge 12$. By construction of ${\cal S}_N^{\rm (n)}$, $a_{i'}$ is such that  $a_{i'} \cap {\cal L} = a_{i'} \cap {\cal S}_N$ if $2 \le i' \le 5$ and $[A_1 + {e}_2, A_2] \cap {\cal L} = [A_1 + {e}_2, A_2] \cap {\cal S}_N$ if $i' = 1$. We can then generate a different ground state in the following way: 

Let $k \in \N$ such that $|a_i| + k \le |a_{i'}| - k - 1$. Detach a portion ${\cal U}$ from ${\cal S}_N$ (as sketched in Figure \ref{fig:nolossii}) consisting of all the atoms in ${\cal S}_N$ which can be connected by at most $k-1$ lattice edges of ${\cal L}$ to some lattice point on $a_i$. Now attach ${\cal U}$---suitably translated and rotated---to $a_{i'}$ (resp.\ $[A_1 + {e}_2, A_2]$ if $i' = 1$) restoring all atomic bonds. Then also the new configuration $\tilde{{\cal S}}_N$ is a ground state. 

Note that the union $\tilde{\cal C}_i \cup \tilde{\cal C}_{i + 1}$ of the corresponding connected components of $\tilde{\cal T} \setminus \tilde{{\cal S}}_N$ (cf.\ Lemma \ref{lemma:structureii}) contains at least $(|a_{i'}| - k - |a_i|) k$ lattice sites. Since $|a_{i'}| - |a_i| \ge 12$, we may assume that 
$$ \frac{1}{4} (|a_{i'}| - |a_i|) \le k \le \frac{1}{3} (|a_{i'}| - |a_i|), $$ 
so that  
$$ \# ( \tilde{\cal C}_i \cup \tilde{\cal C}_{i + 1} ) \ge \frac{1}{6} (|a_{i'}| - |a_i|)^2. $$
On the other hand, by Lemma \ref{lemma:structureii}(iv) and Theorem \ref{theo:HeitmannRadin}(iii) we know that $\# ( \tilde{\cal C}_i \cup \tilde{\cal C}_{i + 1} ) \le C \sqrt{N}$. It follows that $|a_{i'}| - |a_i| \le C N^{\frac{1}{4}}$ (which of course also holds true if $|a_{i'}| - |a_i| \le 12$). Now using that $a_1 \cup \ldots \cup a_6$ is a closed hexagon, we find that indeed for all $i, i' \in \{1, \ldots, 6\}$
$$ \left| |a_{i'}| - |a_i| \right| \le C N^{\frac{1}{4}}. $$

Theorem \ref{theo:HeitmannRadin}(iii) now implies that 
$$ |a_i| = \sqrt{\frac{N}{3}} + {\cal O}(N^{\frac{1}{4}}),~ i = 1, \ldots, 6. $$
Rescaling by $N^{-\frac{1}{2}}$ concludes the proof. \eop

\begin{figure}[h!]
\begin{center}
\setlength{\unitlength}{0.5cm}
\begin{picture}(11,8)
\multiput(3,0)(1,0){5}{\circle*{0.30}}
\multiput(2.5,0.866)(1,0){6}{\circle*{0.30}}
\multiput(2,1.732)(1,0){7}{\circle*{0.30}}
\multiput(1.5,2.598)(1,0){8}{\circle*{0.30}}
\multiput(1,3.464)(1,0){9}{\circle*{0.30}}
\multiput(1.5,4.33)(1,0){8}{\circle*{0.30}}
\multiput(2,5.196)(1,0){7}{\circle*{0.30}}
\multiput(2.5,6.062)(1,0){6}{\circle*{0.30}}
\multiput(3,6.928)(1,0){5}{\circle*{0.30}}

\multiput(8,0)(0.5,0.866){5}{\circle{0.30}}
\multiput(10,3.464)(-0.5,0.866){6}{\circle{0.30}}
\multiput(6.5,7.794)(-1,0){2}{\circle{0.30}}
\end{picture}
\end{center}
\caption{\label{F:hex} Natural candidate for a discrete ground state. Black atoms belong to ${\cal S}_N' \cap \operatorname{conv}(B_1^{(k)}, \ldots, B_6^{(k)})$, white atoms belong to ${\cal S}_N' \setminus \operatorname{conv}(B_1^{(k)}, \ldots, B_6^{(k)})$.} 
\end{figure}
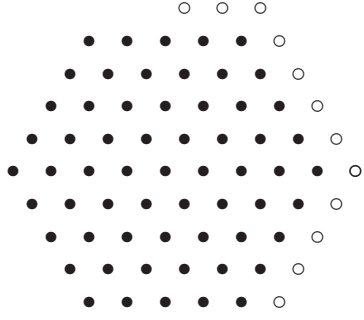

In order to prove by example that the bound obtained in Theorem \ref{theo:rateupper} is sharp, we start by noting that the following natural candidate for an energy minimizing configuration is indeed a ground state. For $k \in \N$ consider the regular hexagon spanned by $B_i = B_i^{(k)}$, where $B_1 = k({e}_1 - {e}_2)$, $B_2 = k {e}_1$, $B_3 = k {e}_2$ and $B_4 = - B_1$, $B_5 = - B_2$, $B_6 = - B_3$. Choose $k \in \N$ such that 
$$ 3 k (k + 1) + 1 \le N < 3 (k + 1)(k + 2) + 1. $$
Now let ${\cal S}_N' \subset \operatorname{conv}(B_1^{(k + 1)}, \ldots, B_6^{(k + 1)}) \cap {\cal L}$ be the configuration with $N$ atoms that satisfies 
$$ {\cal S}_N' \cap \operatorname{conv}(B_1^{(k)}, \ldots, B_6^{(k)}) 
   = {\cal L} \cap \operatorname{conv}(B_1^{(k)}, \ldots, B_6^{(k)}) $$
and, in case ${\cal S}_N' \setminus \operatorname{conv}(B_1^{(k)}, \ldots, B_6^{(k)}) \ne \emptyset$, $B_1^{(k)} + {e}_1 \in {\cal S}_N'$, $B_1^{(k)} + {e}_1 - {e}_2 \notin {\cal S}_N'$ and the set of atoms ${\cal S}_N' \setminus \operatorname{conv}(B_1^{(k)}, \ldots, B_6^{(k)})$ is connected (see Figure \ref{F:hex}). 

The following Proposition is proved in \cite{HeitmannRadin:80}. 

\begin{proposition}[\cite{HeitmannRadin:80}]\label{prop:naturalcandidate}
${\cal S}_N'$ is a ground state.
\end{proposition}

Note that in fact by recourse to Theorem \ref{theo:HeitmannRadin}(i) the proof of this result simply amounts to counting the number of atomic bonds in ${\cal S}_N'$. 


We will now show that the estimate on the rate of convergence in Theorem \ref{theo:rateupper} is optimal.

\begin{theorem} \label{theo:lowerbound}
There exists a constant $c > 0$, independent of $N$, such that for infinitely many $N \in \N$ there are two ground states ${\cal S}_N'$ and $\tilde{{\cal S}}_N$ with $N$ atoms such that 
$$ \inf \left\{ \| \mu_N' - \tilde{\mu}_N(R \cdot + a) \| : R \in O(2), a \in \R^2 \right\} 
   \ge c N^{-\frac{1}{4}}, $$
where $\mu_N' = \frac{1}{N} \sum_{x \in {\cal S}_N'} \delta_{\frac{x}{\sqrt{N}}}$, $\tilde{\mu}_N = \frac{1}{N} \sum_{x \in \tilde{{\cal S}}_N} \delta_{\frac{x}{\sqrt{N}}}$.  
\end{theorem}

\Proof Let ${\cal S}_N'$ be the ground state constructed above (cf.\ Proposition \ref{prop:naturalcandidate}) and suppose $N = 3 k (k + 1) + 2$ for some $k \in \N$. Then 
$$ {\cal S}_N' = ({\cal L} \cap \operatorname{conv}(B_1, \ldots, B_6)) \cup \{ B_1 + {e}_1 \} $$ 
with $B_i = B_i^{(k)}$ as above. Note that for $H : = \overline{h_k} = \operatorname{conv}(B_1, \ldots, B_6)$ we have 
$$ \left\| \mu_N' - \frac{2}{\sqrt{3}} \chi_{N^{-1/2} H} \right\| \le C N^{- \frac{1}{2}}. $$
Now let $m := \left\lfloor \sqrt{\frac{k}{2}} \right\rfloor$ and move the atoms at positions 
$$ B_2 - s {e}_2 + t ({e}_2 - {e}_1), \quad s, t \in \{0, \ldots, m-1\} $$
successively to the points $B_1 + {e}_1 + i {e}_2$, $i = 1, 2, \ldots, m^2$ in such a way that the net loss of atomic bonds is zero. 

\begin{figure}[h!]
\begin{center}
\setlength{\unitlength}{0.8cm}
\begin{picture}(8,5.8)
\put(0,-0.5){%
\setlength{\unitlength}{0.8cm}
\begin{picture}(7,6.2)

\multiput(2,1)(0.01,0){301}{\circle*{0.01}}
\multiput(2,1)(-0.01,0.0173){151}{\circle*{0.01}}
\multiput(0.5,3.598)(0.01,0.0173){151}{\circle*{0.01}}
\multiput(2,6.196)(0.01,0){301}{\circle*{0.01}}
\multiput(5,1)(0.01,0.0173){131}{\circle*{0.01}}
\multiput(5,6.196)(0.01,-0.0173){101}{\circle*{0.01}}
\multiput(6,4.464)(-0.01,-0.0173){21}{\circle*{0.01}}
\multiput(4.6,6.196)(0.01,-0.0173){171}{\circle*{0.01}}

\multiput(5,1)(-0.01,-0.0173){6}{\circle*{0.01}}
\multiput(4.9,0.827)(-0.01,-0.0173){6}{\circle*{0.01}}

\multiput(4.8,0.654)(-0.01,0){12}{\circle*{0.01}}
\multiput(4.6,0.654)(-0.01,0){12}{\circle*{0.01}}
\multiput(4.4,0.654)(-0.01,0){12}{\circle*{0.01}}
\multiput(4.2,0.654)(-0.01,0){12}{\circle*{0.01}}
\multiput(4.0,0.654)(-0.01,0){12}{\circle*{0.01}}
\multiput(3.8,0.654)(-0.01,0){12}{\circle*{0.01}}
\multiput(3.6,0.654)(-0.01,0){12}{\circle*{0.01}}
\multiput(3.4,0.654)(-0.01,0){12}{\circle*{0.01}}
\multiput(3.2,0.654)(-0.01,0){12}{\circle*{0.01}}
\multiput(3,0.654)(-0.01,0){12}{\circle*{0.01}}

\multiput(2.8,0.654)(-0.01,0.0173){6}{\circle*{0.01}}
\multiput(2.7,0.827)(-0.01,0.0173){6}{\circle*{0.01}}

\put(6.4,3.2){\footnotesize $B_2$}
\put(5.15,1){\footnotesize $B_1$}

\qbezier(5.85,5.15)(11,1.7)(5,0.8)
\put(5,0.8){\vector(-4,-1){.01}}
\end{picture}
}
\end{picture}
\end{center}
\caption{\label{fig:noloss} Moving a large set of atoms keeping the number of bonds fixed.}
\end{figure}
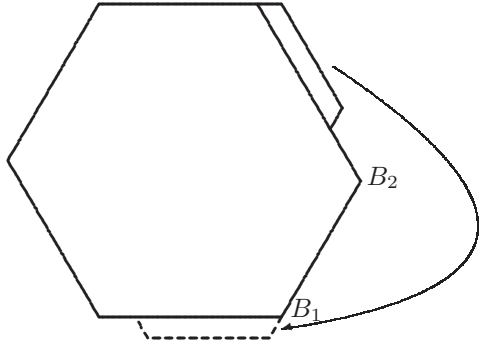

Denoting the new configuration and the corresponding rescaled empirical measure by $\tilde{{\cal S}}_N'$ and $\tilde{\mu}_N'$, respectively, we thus obtain that 
$$ \left\| \tilde{\mu}_N - \tilde{\mu}_N' \right\| \le C N^{-\frac{1}{2}}. $$
Now define $\tilde{{\cal S}}_N$ by detaching the whole set of atoms with positions in 
$$ \{ B_2 - s {e}_2 + t({e}_2 - {e}_1) : s \in \{0, \ldots, m - 1\}, t \in \N \} $$
from $\tilde{{\cal S}}_N'$ and attaching it---after a suitable translation and rotation---to the atoms on the segment $[B_6, B_1]$ (cf.\ Figure \ref{fig:noloss}). 

This can be done such that the number of atomic bonds is preserved, so that we arrive at yet another ground state with corresponding rescaled empirical measure $\tilde{\mu}_N$. 

If $\tilde{H}$ denotes the hexagonal set 
$$ \tilde{H} := \operatorname{conv}(B_1 - m {e}_2, B_2 - m {e}_2, B_3 - m {e}_1, B_4, B_5, B_6 + m ({e}_1 - {e}_2)), $$
then 
$$ \left\| \tilde{\mu}_N - \frac{2}{\sqrt{3}} \chi_{N^{-1/2} \tilde{H}} \right\| \le C N^{- \frac{1}{2}}. $$
Note that each boundary segment of $H$ and $\tilde{H}$ is of length $\sqrt{\frac{N}{3}} + \mathcal{O}(m)$, $m \sim N^{\frac{1}{4}}$. However, while all the segments of $\partial H$ are of equal length, those of $\partial \tilde{H}$ differ in length by an amount $\ge c m \ge c N^{\frac{1}{4}}$. It is not hard to see that indeed 
$$ \inf \{ |H \triangle (R H + a)| : R \in O(2), a \in \R^2 \} \ge c N^{\frac{3}{4}}.$$
(Note e.g.\ that $\operatorname{diam}(\tilde{H}) - \operatorname{diam}(H) \ge \frac{m}{2}$.) The assertion of the theorem now follows by rescaling. \eop 

Note that this degeneracy of the ground state energy is in sharp contrast to the situation with one atom less:

\begin{proposition}\label{prop:unique-gs} Suppose $N = 3 k (k + 1) + 1$ for some $k$. The ground state is (up to rigid motions) uniquely given by 
$$ {\cal S}_N = {\cal L} \cap \operatorname{conv}(B_1^{(k)}, \ldots, B_6^{(k)}). $$
\end{proposition}
So as $N$ grows, the situation changes rapidly between highly degenerate and non-degenerate energy levels. \smallskip 

\Proof Suppose $\tilde{{\cal S}}_N \subset {\cal L}$ is another ground state with $N$ atoms and denote the corresponding set of atomic bonds by $\tilde{\cal B}$, its boundary polygon by $\tilde{P}$. Let $F$ and $\tilde{F}$ be the interiors of $P$ and $\tilde{P}$, respectively. The configurations ${\cal S}_N$ and $\tilde{{\cal S}}_N$ introduce a triangulation of $F$ and $\tilde{F}$, respectively, into $n_F$, respectively, $n_{\tilde{F}}$, triangles of volume $\frac{\sqrt{3}}{4}$. Since $\#{\cal B} = \#\tilde{\cal B}$ and, by Theorem \ref{theo:HeitmannRadin}, $F$ and $\tilde{F}$ are simply connected, we have by Euler's formula $n_F = n_{\tilde{F}}$ and hence $|F| = |\tilde{F}|$.  
Again by Theorem \ref{theo:HeitmannRadin} we have
$$ 
      {\cal H}^1(P) = {\cal H}^1(\tilde{P}).
$$

We now use the following three facts. First, by inspection the above perimeters can be expressed in terms of the following surface energy functional derived in \cite{AFS} as a Gamma-limit of the Heitmann-Radin energy:
$$
   2 {\cal H}^1(P) = \int_{\partial F} \Gamma(\nu), \;\;\; 2 {\cal H}^1(\tilde{P}) = \int_{\partial \tilde{F}} \Gamma(\nu),
$$
where $\nu$ denotes the outward unit normal to the domain of integration and $\Gamma$ is the $\pi/3$-periodic function defined by 
\begin{align}\label{hexenergy}
  \Gamma(\nu) = 2(\nu_2-\nu_1/\sqrt{3}) \mbox{ for }\nu = {-\sin\varphi \choose \cos\varphi}, \varphi\in[0,\pi/3].
\end{align}
Second, such surface energy functionals have a unique minimizer given by the Wulff construction, as proved by 
Taylor (in the language of geometric measure theory)
and Fonseca-M\"uller (in the present language of boundary integrals):
\begin{theorem}\label{theo:TaylorFonsecaMueller} {\rm \cite{Taylor2, FonsecaMueller}} A functional of form
$$
     I(E) = \int_{\partial^\ast E} \Gamma(\nu(x))\, d{\cal H}^{n-1}(x),
$$
with $\Gamma \, : \, S^{n-1}\to[0,\infty)$ continuous and bounded away from zero,
is minimized over sets $E\subset\R^n$ of finite perimeter and volume $c > 0$ if and only if $E$ agrees,
up to translation and up to a set of measure zero, with $\lambda_c W_\Gamma$, where $W_\Gamma$ is
the Wulff set
$$
      W_\Gamma := \Bigl\{x\in\R^n : x\cdot \nu\le \Gamma(\nu) \mbox{ for all }\nu\in S^{n-1}\Bigr\}
$$
and $\lambda_c>0$ is the unique normalization constant such that $\lambda W_\Gamma$ has volume $c$.
\end{theorem}
Finally, for the energy (\ref{hexenergy}), an elementary calculation shows that the
Wulff set is given by the intersection of the six half-spaces $x\cdot\nu \le \Gamma(\nu)$ for
the minimizing normals $\nu(\frac{2\pi j}{6})$, $j=0,\ldots,5$. 

Since $\overline{F} = \lambda_{|F|} W_{\Gamma}$, this establishes that both $F$ and $\tilde{F}$ are minimizers of $I$ subject to the same fixed volume, so that by 
Theorem \ref{theo:TaylorFonsecaMueller} they coincide and therefore 
${\cal S}_N = \tilde{{\cal S}}_N$ (up to translation). \eop

%
%

\section*{Acknowledgment}
I am very grateful to Yuen AuYeung and Gero Friesecke for all our stimulating discussions on discrete particle systems and emerging Wulff shapes. The results of this paper and their exposition have benefitted a lot from their input and help in preparing this manuscript.

%
%

\bibliographystyle{alpha}

\end{document}